\documentclass[12pt,reqno]{amsart} 
\usepackage{amsmath,amssymb,amsthm,mathrsfs}
\usepackage[english]{babel}
\usepackage[utf8]{inputenc}
\usepackage{xcolor}
\newtheorem{theorem}{Theorem}[section]
\newtheorem{proposition}[theorem]{Proposition}
\newtheorem{corollary}[theorem]{Corollary}

\theoremstyle{definition}
\newtheorem{definition}[theorem]{Definition}
\newtheorem{example}[theorem]{Example}
\newtheorem{remark}[theorem]{Remark}

\usepackage{multirow} 
\usepackage{color}

\usepackage{graphicx}
\usepackage[left=3cm,right=3cm]{geometry}
\usepackage{bigints}
\usepackage{comment}

\def\r{\mathbb R}
\def\h{\mathbb H}
\def\s{\mathbb S}
\def\sl{\mathrm{SL}(2,\r)} 
 \def\nn{\mathcal N}
  \def\aa{\mathcal A}
   \def\kk{\mathcal K}
   
\begin{document}

\title[Translators of the mean curvature flow in $\mathrm{SL}(2,\mathbb{R})$]{Translators of the mean curvature flow in the special linear group   $\sl$}
\author{Rafael L\'opez, Marian Ioan Munteanu}
\address{ Departamento de Geometr\'{\i}a y Topolog\'{\i}a\\  Universidad de Granada. 18071 Granada, Spain}
\email{rcamino@ugr.es}
\address{University 'Al. I. Cuza' of Iasi, Faculty of Mathematics, Bd. Carol I, no. 11, 700506 Iasi, Romania}
\email{marian.ioan.munteanu@gmail.com}

\begin{abstract}
Translators in the special linear group $\sl$ are surfaces whose mean curvature $H$ and unit normal vector $N$ satisfy $H=\langle N,X\rangle$, where $X$ is a fixed Killing vector field. In this paper we study and classify those translators that are invariant by a one-parameter group of isometries. By  the Iwasawa decomposition, there are three types of such groups. The dimension of the Killing vector fields is $4$ and an   exhaustive discussion is done for each one of the Killing vector fields and each of the invariant surfaces. In some cases, explicit parametrizations of translators are obtained.   
 \end{abstract}


\maketitle

\section{Introduction and preliminaries}

In the theory of mean curvature flow (MCF for short)  in Euclidean space $\r^3$, translators of the MCF are  
 surfaces that evolve purely by translations of $\r^3$. A translator  $\Sigma\subset\r^3$ is characterized by
\begin{equation}\label{eq1}
H=\langle N,\textbf{v}\rangle,
\end{equation}
where $H$ and $N$ is the mean curvature and unit normal of $\Sigma$ and $\textbf{v}$ is the direction of translations of $\r^3$. The role of translators is important because they are, after rescaling,  a type of singularities of the MCF  according to Huisken and Sinestrari  \cite{hs}.    The simplest example is any plane  parallel to $\textbf{v}$. Other examples of translators are those that  are invariant by a one-parameter group of  isometries of $\r^3$. If the translator is invariant in one   spatial direction of $\r^3$, the translator is the grim reaper. On the other hand, if the translator is a surface of revolution, then the rotation axis is parallel to $\textbf{v}$. In such a case there are two types of translators depending whether the surface intersects the rotation axis (bowl soliton) or not (wing-like rotational translators).   Recently, there is a great interest in extending   the notion of translators and solitons in general  of the MCF  in other homogenous spaces.    Without to be complete, we refer: hyperbolic space \cite{bl1,bl2,li,mr}; the product $\h^2\times\r$ \cite{bue1,bue2,bl3,lipi};  the product $\s^2\times\r$ \cite{lm1}; the Sol space \cite{pi1}; and the Heisenberg group \cite{pi2}.

In this paper we   study    translators of the MCF  in  the special linear group $\sl$. This space is viewed as a   homogeneous space  equipped with a  canonical left-invariant Riemannian metric  whose group of isometries is of dimension $4$. For the definition of translator in $\sl$, we extend the notion in \eqref{eq1}   replacing $\textbf{v}$ by a Killing vector field of $\sl$.

\begin{definition}\label{def1} Let $X\in\mathfrak{X}(\sl)$ be a Killing vector field. A surface $\Sigma$ in $\sl$ is said to be a {\it $X$-translator}  if its mean curvature $H$ and unit normal vector $N$ satisfy
\begin{equation}\label{eq2}
H=\langle N,X\rangle,
\end{equation}
\end{definition}

Let 
$$\sl=\left \{ \left (
\begin{array}{cc}
a & b \\ c & d
\end{array}
\right ): a,b,c,d \in \r, ad-bc=1\right \}.$$
For the definition of the Riemannian metric, we recall   the so-called    Iwasawa decomposition \cite{iw}. In the group $\sl$ there are the following three (one-dimensional) subgroups:

\begin{align*} \nn&=\left \{
\left (\begin{array}{cc}
1 & x \\ 0 & 1 
\end{array}
\right ): x \in \r\right \},  
 \\
\aa&=\left \{\left (
\begin{array}{cc}
\sqrt{y} & 0 \\ 0 & 1/\sqrt{y}
\end{array}
\right ): y\in\r_+
\right \},  \\
\kk&=\left \{\left (\begin{array}{cc}
\cos \theta & \sin \theta \\
-\sin \theta & \cos \theta 
\end{array}\right ): 
  \theta\in\r 
\right \} .
\end{align*}
Notice the isomorphisms $\nn= (\r,+)$, $\aa=(\r_+,\cdot)$ and $\kk=(\s^1,\cdot)$. By the Iwasawa decomposition, denoted by NAK decomposition,   for every $A\in\sl$ there is a unique representation of $A$ given by $A=nak$, where $n\in\mathcal{N}$, $a\in\aa$ and $k\in\mathcal{K}$.  This   allows to give  global coordinates  $(x,y,\theta)$ in $\sl$ by means of
\begin{equation} 
\label{coord}
(x,y,\theta)\in\r^3\longmapsto
\left (\begin{array}{cc}
1 & x \\ 0 & 1 
\end{array}
\right )
\left (\begin{array}{cc}
\sqrt{y} & 0 \\ 0 & 1/\sqrt{y}
\end{array}\right )
\left (\begin{array}{cc}
\cos \theta & \sin \theta \\
-\sin \theta & \cos \theta 
\end{array}\right )\in \sl.
\end{equation}
With respect to these coordinates, let $\{\partial_x,\partial_y,\partial_\theta\}$ be the canonical basis of  $\mathfrak{X}(\sl)$.

\begin{remark}\label{rm1} If $A\in\sl$, then the characteristic polynomial of $A$ is 
$\lambda^2-\mbox{trace}(A)\lambda+1$. This     distinguishes the matrices of $\sl$ in three types depending on the number of roots of this polynomial. This classification is equivalent to the NAK decomposition. Indeed,  if $|\mbox{trace}(A)|=2$, there is a unique double eigenvalue; if $|\mbox{trace}(A)|>2$, there are two distinct real eigenvalues and if   $|\mbox{trace}(A)|<2$, there are no real eigenvalues. Examples of such matrices are, respectively, that of  $\nn$, $\aa $ and $\kk$. In the literature,  the elements of the subgroups $\nn$, 
$\aa$ and $\kk$ are also called parabolic, hyperbolic and elliptic matrices, respectively \cite{ki}. 
\end{remark}

Let $\h^2(-4)$ be the hyperbolic plane of constant curvature $-4$ and its   upper half plane model  
 \[
\mathbb{H}^2(-4)=\left(\r^2_+, \frac{dx^2+dy^2}{4y^2}
\right).
\]
The special linear group $\sl$ 
acts transitively and isometrically on $\h^2(-4)$ by the 
linear fractional transformation  
\[
\left(A=\left(
\begin{array}{cc}
a & b\\
c & d
\end{array}
\right),(x,y)=x+iy=z\right)\mapsto A\cdot z=\frac{az+b}{cz+d}.
\]
 The isotropy subgroup of $\sl$ at $i=(0,1)$ is the subgroup $\kk$. In the terms of the NAK decomposition, the natural projection 
\begin{equation}\label{pro}
\pi:\sl
\to \sl/\kk=\mathbb{H}^2(-4),\quad \pi(x,y,\theta)=(x,y).
\end{equation}
 The mapping 
$$\psi:\h^2(-4) \times 
\s^1 \rightarrow \sl,$$
$$\psi(x,y,\theta)=
\begin{pmatrix} 1 & x \\ 0 & 1 
\end{pmatrix}
\begin{pmatrix}
\sqrt{y} & 0 \\ 0 & 1/\sqrt{y} 
\end{pmatrix} 
\begin{pmatrix}
\cos \theta & \sin \theta \\ 
-\sin \theta & \cos \theta
\end{pmatrix}
$$
is a diffeomorphism. In particular,  the space $\sl$ is topologically   the  open solid torus $\mathbb{D} \times \s^1$.  If we endow $\sl$ by the metric $\langle~,~\rangle$ which makes $\psi$ an isometry,  the expression of $\langle,\rangle$ is    
$$\langle ~, ~\rangle=\frac{dx^2+dy^2}{4y^2}+\left(d\theta+\frac{dx}{2y}\right)^2.$$
With this metric, the   projection $\pi$ defined in \eqref{pro} becomes a Riemannian submersion.

The space of Killing vector fields in $\sl$ is of dimension $4$ and it is   generated by 
\begin{equation}\label{kk}
\{\partial_x,  \partial_\theta, 
x\partial_x+y\partial_y,   \frac12(x^2-y^2)\partial_x+xy\partial_y\}.
\end{equation}

In order to give examples of translators of $\sl$, and following the motivation from  the Euclidean case,  we will assume that the surface is invariant by a one-parameter group of isometries.   The NAK decomposition allows us to give the following definitions.

 \begin{definition} 
 Let $\Sigma$ be an immersed  surface in $\sl$. We say that $\Sigma$ is 
 $\nn$-invariant (resp.  $\aa$-invariant, $\kk$-invariant) if $\Sigma$ 
 is invariant under the left translations of the subgroup $\nn$ 
 (resp.  $\aa$, $\kk$). Also,   $\kk$-invariant  surfaces are called   rotational surfaces. 
 \end{definition}

Invariant surfaces  in $\sl$ with constant mean curvature or constant Gauss curvature  have been studied in  \cite{er,in,ko,mu,to}. 
 
Once we have established the  definition of invariant surface, the work ahead is the classification of invariant $X$-translators depending of the Killing vector field $X$ of \eqref{kk}. The paper is organized in sections according the Killing vector field $X$. In each section, namely, Sects. \ref{s3}, \ref{s4}, \ref{s5} and \ref{s6}, we will study and classify the $X$-translators that are invariant by each of the three subgroups $\nn$, $\aa$ and $\kk$. Previously, in Sect. \ref{s2}, we compute the unit normal vector $N$ and the mean curvature $H$ of the invariant surfaces. These computations are needed to study the translator equation \eqref{eq2}.

 By the variety of vector fields and invariant surfaces, we summarize in Table \ref{table1} the results of classification obtained in this paper. In the table, by explicit parametrization, we mean that we obtain a parametrization of the surface by known functions. Other surfaces that we obtain are    those where one of the coordinates $x$, $y$ or $\theta$ in the Iwasawa decomposition are constant. Let $\Sigma_{x_0}$, $\Sigma_{y_0}$ and $\Sigma_{\theta_0}$ be the corresponding surfaces, respectively. Finally, in the case of $\kk$-invariant translator, by ODE we mean that we obtain the differential equation that describes the generating curve of the surface. In general, this equation is difficult to study in all its generality.

\begin{table}[ht]
\centering
\begin{tabular}{|c|c|c|c|}
\hline
&$\nn$-surfaces&$\aa$-surfaces&$\kk$-surfaces\\
\hline
$\partial_x$& \mbox{explicit parametrization}&$\Sigma_{\theta_0}$ &\mbox{description} \\\hline
$\partial_\theta$&\mbox{$\Sigma_{\theta_0}$, explicit parametrization} & $\Sigma_{x_0}$ &\mbox{minimal surface} \\\hline
$x\partial_x+y\partial_y$& \mbox{explicit parametrization} & $\Sigma_{x_0=0}$, $\Sigma_{\theta_0}$ & ODE\\\hline
$ \frac12(x^2-y^2)\partial_x+xy\partial_y$& $\Sigma_{\theta_0}$ & $\Sigma_{\theta_0}$ & ODE\\\hline
\end{tabular}
\caption{Classification of the invariant $X$-translators. }\label{table1}
\end{table}
We obtain two direct consequences.

\begin{corollary} All   $\nn$-invariant translators have   explicit parametrizations by known functions.   
\end{corollary}

\begin{corollary} The only   $\aa$-invariant translators are of type  $\Sigma_{x_0}$ or $\Sigma_{\theta_0}$ depending on the case.
\end{corollary}

 \section{The  mean curvature of invariant surfaces}\label{s2}

In this section, we compute the unit normal vector $N$ and the mean curvature $H$ of surfaces invariant by each one of the three subgroups $\nn$, $\aa$ and $\kk$ of $\sl$. Part of the computations of this section have appeared in \cite{ko,mu}. By completeness of the paper and for the subsequent study of the $X$-translators, we recall them.   Consider in $\sl$   the orthonormal frame  $B=\{e_1,e_2,e_3\}$ defined by  
\begin{equation}\label{framefield}
e_1=2y\partial_x
-\partial_\theta,\ e_2=2y\partial_y,\  e_3=\partial_\theta.
\end{equation}

The Levi-Civita connection $\nabla$ of the metric $\langle~,~\rangle$ of $\sl$ is given by the relations
$$\begin{array}{lll}
\nabla_{e_1}e_1=2e_2,& \nabla_{e_1}e_2=-2e_1-e_3,&\nabla_{e_1}e_3=e_2,\\
\nabla_{e_2}e_1=e_3,&\nabla_{e_2}e_2=0,&\nabla_{e_2}e_3=-e_1,\\
\nabla_{e_3}e_1=e_2,&\nabla_{e_3}e_2=- e_1,&\nabla_{e_3}e_3=0.
\end{array}
$$
The expressions of Killing vector fields \eqref{kk} in terms of the above basis are  
\begin{equation*}
\begin{split}
\partial_x&=\frac{1}{2y}(e_1+e_3),\\
\partial_\theta&=e_3,\\
x\partial_x+y\partial_y& =\frac{1}{2y}(xe_1+y e_2+xe_3),\\
\frac12(x^2-y^2)\partial_x+xy\partial_y&=\frac{1}{2y}\left(\frac12(x^2-y^2) e_1+xy e_2+\frac12(x^2-y^2)e_3\right).
\end{split}
\end{equation*}

We now give parametrizations of the invariant surfaces in order to compute $N$ and $H$. To have a consistent notation, if $\Psi=\Psi(s,t)$ are local coordinates on the invariant surface, the parameter $s$ will be assigned for the generating curve of the surface, whereas   $t$ will denote the parameter of the group.

\begin{proposition}
\label{prop-2.1} 
An  $\nn$-invariant surface of $\sl$ can be parametrized by 
\begin{equation}
\label{p-Nsurf1}
\Psi(s,t)=\left (\begin{array}{cc}
1 & t \\ 0 & 1 
\end{array}
\right )
\left (\begin{array}{cc}
\sqrt{y(s)} & 0 \\ 0 & 1/\sqrt{y(s)}
\end{array}\right )
\left (\begin{array}{cc}
\cos \theta(s) & \sin \theta(s) \\
-\sin \theta(s) & \cos \theta(s) 
\end{array}\right ),
\end{equation}
where $t\in \r$, $s\in I\subset \r$.  
The generating curve is $\alpha(s)=(y(s),\theta(s))$. 
The unit normal $N$ and the mean curvature $H$ are 
\begin{equation}\label{nh-n}
\begin{split}
N&=\frac{y'}{\sqrt{2}\Phi}(e_1-e_3)+\frac{ \sqrt{2}y\theta'}{\Phi}e_2,\\
H&=\frac{\sqrt{2}y^2}{\Phi^3}\left(\theta'y''-y'\theta''+2y\theta'^3\right),
\end{split}
\end{equation}
where $  \Phi=\sqrt{y'^2+2y^2\theta'^2}$. 
If we take    $\alpha(s)$ such that 
\begin{equation}\label{nh-0}
\begin{split}
y'(s)&=\sqrt{2}y(s)\cos\varphi(s),\\
\theta'(s)&=\sin\varphi(s),
\end{split}
\end{equation}
then 
\begin{equation}\label{nh-n2}
\begin{split}
N&=\frac{1}{\sqrt{2}}\cos\varphi(e_1-e_3)+\sin\varphi e_2,\\
 H&=-\frac{\varphi'}{\sqrt{2}}+\sin\varphi.
 \end{split}
 \end{equation}

\end{proposition}

\begin{proposition} An $\aa$-invariant surface   of $\sl$ can be parametrized by 
\begin{equation}
\label{p-conoid}
\Psi(s,t)=\left (\begin{array}{cc}
1 & x(s) \\ 0 & 1 
\end{array}
\right )
\left (\begin{array}{cc}
\sqrt{t} & 0 \\ 0 & 1/\sqrt{t}
\end{array}\right )
\left (\begin{array}{cc}
\cos \theta(s) & \sin \theta(s) \\
-\sin \theta(s) & \cos \theta(s) 
\end{array}\right ),
\end{equation}
where $t\in \r^+$, $s\in\r$. The generating curve is $\alpha(s)=(x(s),\theta(s))$. The unit normal $N$ and the mean curvature $H$ are 
\begin{equation}\label{nh-a}
\begin{split}
N&=\frac{1}{\Phi}\left(-\left(x'+2t\theta'\right)e_1+x'e_3\right)\\
H&=\frac{2t^2}{\Phi^3}(x'\theta''-\theta'x''),
\end{split}
\end{equation}
where 
$$ \Phi=\sqrt{(x'+2t \theta')^2+x'^2}.\quad
$$
\end{proposition}

\begin{proposition} A rotational surface of $\sl$ can be parametrized by 
\begin{equation}
\label{p-rot}
\Psi(s,t)=\left (\begin{array}{cc}
1 & x(s) \\ 0 & 1 
\end{array}
\right )
\left (\begin{array}{cc}
\sqrt{y(s)} & 0 \\ 0 & 1/\sqrt{y(s)}
\end{array}\right )
\left (\begin{array}{cc}
\cos t & \sin t \\
-\sin t & \cos t 
\end{array}\right ),
\end{equation}
where $s\in I\subset\r$, $t\in\r$. The generating curve is 
$\alpha(s)=(x(s),y(s))\in{\mathbb{H}}^2(-4)$.
The unit normal $N$ and the mean curvature $H$ are 
\begin{equation}\label{nh-k}
\begin{split}
N&=\frac{1}{ \Phi}\left(-y' e_1+x' e_2\right)\\
H&=\frac{1}{\Phi^3} (y(x'y''-x''y')+x'\Phi^2),
\end{split}
\end{equation}
where $ \Phi=\sqrt{x'^2+y'^2}$.  If $\alpha$ is parametrized by arclength, then there is  $\varphi=\varphi(s)$ such that
\begin{equation}\label{xy}
\begin{split}
x'(s)&=2y(s)\cos\varphi(s),\\
 	y'(s)&=2y(s)\sin\varphi(s),
\end{split}
\end{equation}
which implies
\begin{equation}\label{nh-k2}
\begin{split}
N&=-\sin\varphi(s) e_1+\cos\varphi(s) e_2 \\
H&=\frac{\varphi'(s)}{2}+\cos\varphi(s),
\end{split}
\end{equation}
\end{proposition}

We show three particular examples of invariant surfaces which are defined fixing one of the three coordinates $x$, $y$ or $\theta$ in $\sl$.

 \begin{example}  
\label{ex:Sigma-x} 
 Let $\Sigma_{x_0}$ be the   surface in $\sl$ defined by
 \begin{equation}
 \label{S:Sigma-x}
\Psi(s,t)=\left (\begin{array}{cc}
1 & x_0 \\ 0 & 1 
\end{array}
\right )
\left (\begin{array}{cc}
\sqrt{t} & 0 \\ 0 & 1/\sqrt{t}
\end{array}\right )
\left (\begin{array}{cc}
\cos s & \sin s \\
-\sin s & \cos s 
\end{array}\right ).
\end{equation}
The surface $\Sigma_{x_0}$ is a Hopf cylinder over a geodesic in the hyperbolic plane
${\mathbb{H}}^2(-4)$, hence it is both minimal and flat.
The induced metric is $g_{x_0}=ds^2+\frac{dt^2}{4t^2}$.
The surface $\Sigma_{x_0}$ is both $\mathcal{A}$-invariant and $\mathcal{K}$-invariant.
As the unit normal is $N=-e_1$, the surface $\Sigma_{x_0}$ is a   translator with 
respect to $\partial_\theta$ (for any $x_0$) and the surface $\Sigma_{x_0=0}$ is a  translator 
with respect to $x\partial_x+y\partial_y$, too.
 
 \end{example}

 \begin{example}  
\label{ex:Sigma-y} 
Let $\Sigma_{y_0}$ be the surface in $\sl$ defined by 
 \begin{equation}
 \label{S:Sigma-y}
\Psi(s,t)=\left (\begin{array}{cc}
1 & t \\ 0 & 1 
\end{array}
\right )
\left (\begin{array}{cc}
\sqrt{y_0} & 0 \\ 0 & 1/\sqrt{y_0}
\end{array}\right )
\left (\begin{array}{cc}
\cos s & \sin s \\
-\sin s & \cos s 
\end{array}\right ).
\end{equation}
The surface $\Sigma_{y_0}$ is a Hopf cylinder over a Riemannian circle in the 
hyperbolic plane ${\mathbb{H}}^2(-4)$, hence it is flat and of constant mean curvature
$H=1$. The induced metric is $g_{y_0}=\frac{dt^2}{4y_0^2}+(ds+\frac{dt}{2y_0})^2$.
The surface $\Sigma_{y_0}$ is both $\mathcal{N}$-invariant and $\mathcal{K}$-invariant.  The unit normal is $N=e_2$ and $H=1$.  The product of $N$ by the first two vector fields in \eqref{kk} is $0$, and by the last two ones is $\frac12$. Thus  $\Sigma_{y_0}$ is not a translator.
 \end{example}

 \begin{example}  
\label{ex:Sigma-theta} 
 Let $\Sigma_{\theta_0}$ be the     surface in $\sl$ defined by 
 \begin{equation}
 \label{S:Sigma-theta}
\Psi(s,t)=\left (\begin{array}{cc}
1 & s \\ 0 & 1 
\end{array}
\right )
\left (\begin{array}{cc}
\sqrt{t} & 0 \\ 0 & 1/\sqrt{t}
\end{array}\right )
\left (\begin{array}{cc}
\cos \theta_0 & \sin \theta_0 \\
-\sin \theta_0 & \cos \theta_0 
\end{array}\right ).
\end{equation}
The surface $\Sigma_{\theta_0}$ is both $\mathcal{N}$-invariant and 
$\mathcal{A}$-invariant surface. 
The induced metric is $g_{\theta_0}=\frac{2ds^2+dt^2}{4t^2}$ and its curvature
is constant $-4$, that is $\Sigma_{\theta_0}$ is the hyperbolic plane 
$\mathbb{H}^2(-4)$.
Moreover, $\Sigma_{\theta_0}$ is a minimal surface in $\sl$ by \eqref{nh-n}. 
The unit normal is $N=\frac{1}{\sqrt{2}}(e_1-e_3)$. Therefore, $\Sigma_{\theta_0}$ is a 
  translator with respect to  $\partial_\theta$, $x\partial_x+y\partial_y$ and 
$\frac12(x^2-y^2)\partial_x+xy\partial_y$, respectively.
 \end{example}
 
 We end this section with a particular example of rotational surface. 
 
 \begin{example} \label{ex2}
 Let $\Sigma$ be a rotational surface whose generating curve $\alpha$ is a straight-line. Then $\alpha$ is parametrized by \eqref{xy} where the function $\varphi$ is constant.  
 \begin{enumerate}
 \item Case   $\sin\varphi=0$. Without loss of generality, we can suppose $\varphi(s)=0$. An integration of \eqref{xy} gives $\alpha(s)=(2c_2s+c_1,c_2)$, $c_1,c_2\in\r$, $c_2>0$. Then $H=1$ and it is not difficult to check that $\Sigma$ is not a $X$-translator for any vector field $X$ of \eqref{kk}.
 \item Case   $\sin\varphi\not=0$. The solution of \eqref{xy} is $\alpha(s)=c_2 e^{2s\sin\varphi}(\cot\varphi,1)$, $c_1,c_2\in\r$, $c_2>0$. Now $H=\cos\varphi$ and $\Sigma$ is not a $X$-translator with the first and fourth vector fields of \eqref{kk}. For the vector fields $\partial_\theta$ and $x\partial_x+y\partial_y$, we have 
  $\langle N,\partial_\theta\rangle=\langle N,x\partial_x+y\partial_y\rangle=0$. Thus $\Sigma$ is a $\partial_\theta$-translator and  a $x\partial_x+y\partial_y$-translator if and only if $\cos\varphi=0$, that is, $\varphi(s)=\pi/2$. 
 \end{enumerate}
 \end{example}
 
\section{Translators by the vector field $\partial_x$}\label{s3}
Consider the Killing vector field $\partial_x$. We know that with respect to $B$, this vector field   is 
$$\partial_x=\frac{1}{2y}(e_1+e_3).$$

\begin{figure}[hbtp]
\begin{center}
\includegraphics[width=.3\textwidth]{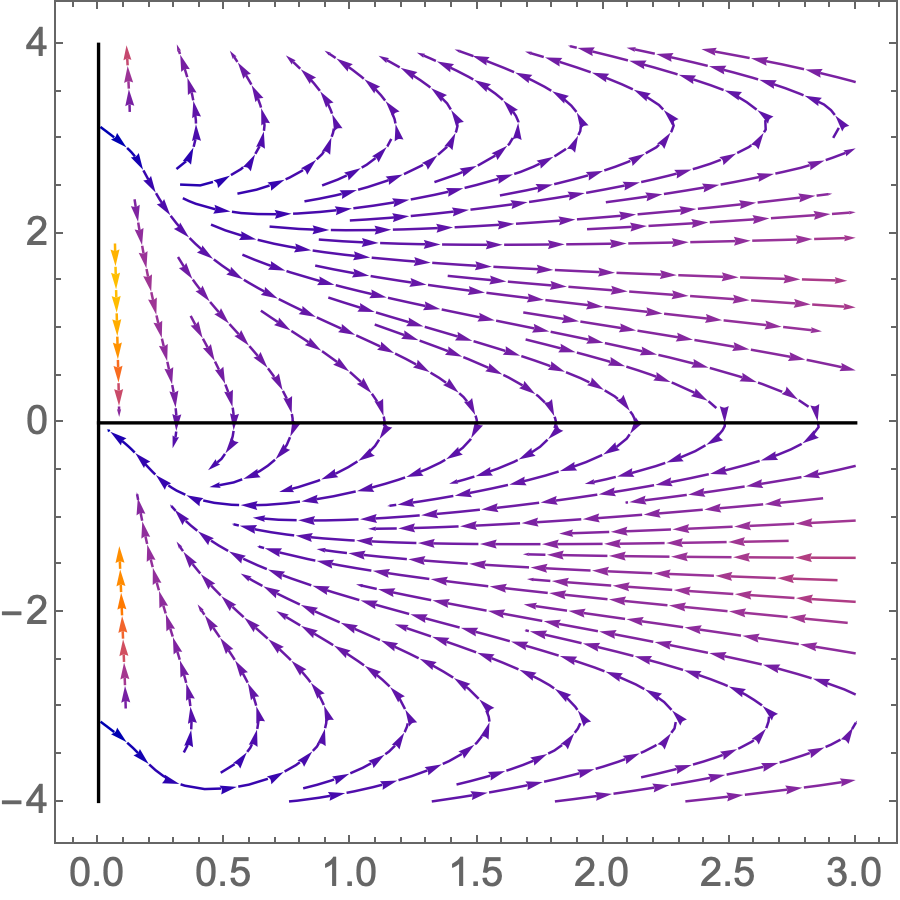}\qquad
\includegraphics[width=.5\textwidth]{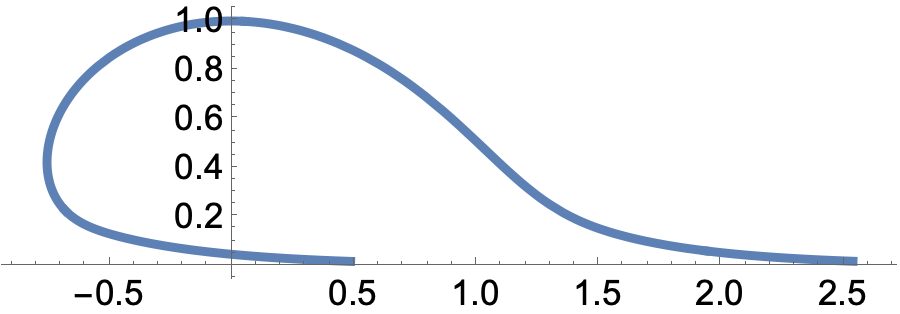}
\end{center}
\caption{Left: the phase portrait of the nonlinear autonomous system \eqref{as}. Right: a solution of \eqref{rt} for initial conditions $\alpha(0)=(0,1)$.}\label{fig1}
\end{figure}

\begin{theorem} 
\begin{enumerate}
\item Let $\Sigma$ be an $\nn$-invariant $\partial_x$-translator. 
Then $\Sigma$ is  either a surface of type $\Sigma_{\theta_0}$ given by \eqref{S:Sigma-theta}, or a minimal surface parametrized by \eqref{p-Nsurf1} with $\theta(s)=s$ and 
$y(s)=c_1\cos(\sqrt{2}s)+c_2\sin(\sqrt{2}s)$, where $c_1$ and $c_2$ are real constants. 
The interval for $s$ is such that $y(s)>0$.
\item The  only $\aa$-invariant  $\partial_x$-translators are the surfaces of type 
$\Sigma_{\theta_0}$.
\item Let $\Sigma$ be a rotational surface whose generating curve is parametrized by \eqref{xy}. If $\Sigma$ is a  $\partial_x$-translator, then  
\begin{equation}\label{rt}
\varphi'=-\frac{\sin\varphi+2y\cos\varphi}{y}. 
\end{equation}
The generating curve $\alpha$ is a bi-graph over the line $y=0$ and converges to it as $s\rightarrow\infty$. See Fig. \ref{fig1}.

 \end{enumerate}
\end{theorem}

\begin{proof} 
\begin{enumerate}

\item  From \eqref{nh-n}, we have $\langle N,\partial_x\rangle=0$, hence the surface is minimal.   If $\theta=\theta_0$ is a constant function $\theta_0$ then $\Sigma=\Sigma_{\theta_0}$ described in Ex. ~\ref{ex:Sigma-theta}. Otherwise, that is, if $\theta'\neq0$ then we can take $\theta(s)=s$. If this is the case, the
equation $H=0$ implies $y''(s)+2y(s)=0$. The solution is $y(s)=c_1\cos(\sqrt{2}s)+c_2\sin(\sqrt{2}s)$, $c_1,c_2\in\r$.


\item 
From \eqref{nh-a}, we have 
$\langle N,\partial_x\rangle=-\dfrac{\theta'}{\Phi}$. Thus an
$\aa$-invariant  $\partial_x$-translator must satisfy
$$
\frac{2t^2}{\Phi^3}(x'\theta''-\theta'x'')=-\frac{\theta'}{\Phi}.
$$
This equation is equivalent to 
\[2 t^2(x'\theta''-\theta' x'')+\theta'x'^2 +\theta' (x'+2t \theta')^2=0.
\]
 Writing this  equation   as a polynomial  equation on the variable $t$, we have
 $$ t^2(x'\theta''-\theta' x''+2\theta'^3)+2x'\theta'^2 t+x'^2\theta'=0.$$
Then all coefficients must vanish. First, $x'^2\theta'=0$. This implies that $\theta'(s)=0$ for all $s$ or $x'(s)=0$ for all $s$. 
However,  $x'$ cannot be identically $0$, otherwise $\theta'(s)$ also vanishes 
for all $s$ and this is not allowed by regularity. Thus  $\theta$ is a constant function, obtaining the surface $\Sigma_{\theta_0}$ of Ex. \ref{ex:Sigma-theta}. 
 

 \item Using  \eqref{nh-k2}, we have  $\langle N,\partial_x\rangle =-\frac{\sin\varphi}{2y}$, obtaining Eq. \eqref{rt}.  Then the generating curve $\alpha(s)=(x(s),y(s))$ is given by the ordinary system formed by the two equations \eqref{xy} together \eqref{rt}. Since the first equation of   \eqref{xy}  can be obtained if we solve the second equation and \eqref{rt}, then it is enough to consider the autonomous system
 \begin{equation}\label{as}
 \left\{\begin{split}
 y'&= 2y\sin\varphi\\
 \varphi'&=-\frac{\sin\varphi(s)}{y}-2\cos\varphi.
 \end{split}\right.
 \end{equation}
 The phase portrait is shown in Fig. \ref{fig1}. The phase plane is $A=\{(y,\varphi):y>0,\varphi\in (-\pi,\pi)\}$. The trajectories cannot intersect by uniqueness of solutions. Each trajectory passing through the point $(y_0,0)$ converges to $(0,0)$ as $s\to\infty$. In particular, the function $y(s)\to 0$ as $s\to\infty$.  This ODE system   appeared in \cite[Th. 5.1]{bl3} in the study of $p$-grim reapers in the space $\h^2\times\r$. We refer there for details.  
 \end{enumerate}
\end{proof}

\section{Translators by the vector field $\partial_\theta$}
\label{s4}
Consider the Killing vector field $\partial_\theta$. This vector field coincides with $e_3$.

\begin{theorem} 
\begin{enumerate}
\item Let $\Sigma$ be an   $\nn$-invariant surface whose generating curve is   parametrized by 
 \eqref{nh-0}. If $\Sigma$ is a $\partial_\theta$-translator, then  either
$\alpha(s)=(c_1e^{2s/\sqrt{3}},-\frac{s}{\sqrt{3}}+c_2)$, where $c_1>0$, 
$c_2\in\mathbb{R}$,  or $\alpha (s)=(y(s),\theta(s))$, 
with
\begin{equation}
\label{nh-1} 
\begin{split}
y(s)&=c_1\exp\left(\frac{2\sqrt{2}}{3}  \Lambda(s)+
	\frac{2}{3}\psi(s) \right)\\
\theta(s)&=\frac{2\sqrt{2}}{3}\Lambda(s)-\frac{1}{3}\psi(s)+c_2,
\end{split}
\end{equation}
 and 
 $ \Lambda(s)=\arctan\left(\tanh\frac{s\sqrt{3}}{2}\right)$, $\psi(s)=\log \cosh(\sqrt{3}),$
 where $c_1>0$, $c_2\in\mathbb{R}$,

\item The  only $\aa$-invariant $\partial_\theta$-translators are the surfaces of type 
$\Sigma_{x_0}$. 
\item The only rotational $\partial_\theta$-translators are minimal surfaces.
\end{enumerate}
\end{theorem}

\begin{proof} 
\begin{enumerate}
\item 
From \eqref{nh-n2}, Eq.   \eqref{eq2} is 
$$\varphi'=\cos\varphi+\sqrt{2}\sin\varphi.$$
A first solution appears when $\varphi$ is a constant function. This occurs when 
$\tan\varphi=-1/\sqrt{2}$. If we consider $\varphi\in[-\frac{\pi}{2},\frac{\pi}{2}]$
it is immediate from \eqref{nh-0} that 
$y(s)=c_1e^{2s/\sqrt{3}}$ and $\theta(s)=-\frac{s}{\sqrt{3}}+c_2$, 
where the constants $c_1$ and $c_2$ are obtained from the initial condition. 
By \eqref{nh-n2} the surface $\Sigma$ is of constant mean curvature $H=-1/\sqrt{3}$.

In case that $\varphi$ is not constant let us fix the initial conditions for $\varphi$
(equivalently to a translation in the parameter $s$ that does not affect the assumption
for the parametrization of $\alpha$)  
$\varphi(0)=\varphi_0=\arctan\sqrt{2}\in(0,\frac{\pi}{2})$. By separation of variables, after integrating, we obtain 
$$
\frac{1+\sin(\varphi-\varphi_0)}{\cos(\varphi-\varphi_0)}=e^{s\sqrt{3}}.
$$
We obtain
$$
\varphi(s)=\arctan\sqrt{2}+
	2\arctan\left(\tanh\frac{s\sqrt{3}}{2}\right).
$$
With this value of $\varphi$, it is possible to integrate \eqref{nh-0} obtaining \eqref{nh-1}. See   Fig. ~\ref{fig:4-1}  the generating curve $\alpha(s)$ in the   $y\theta$-plane.
\begin{figure}[hbtp]
\begin{center}
\includegraphics[width=.4\textwidth]{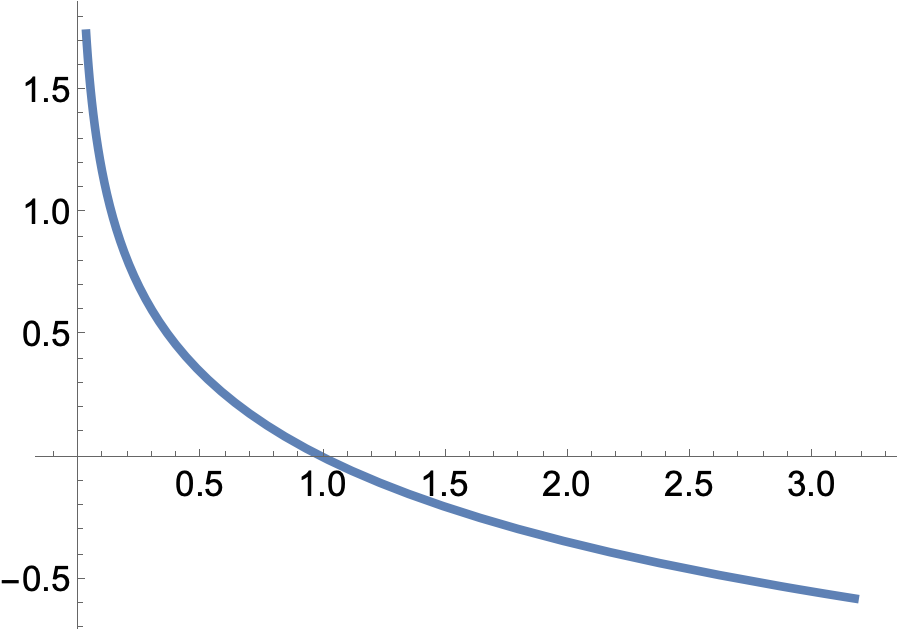}\qquad
\includegraphics[width=.2\textwidth]{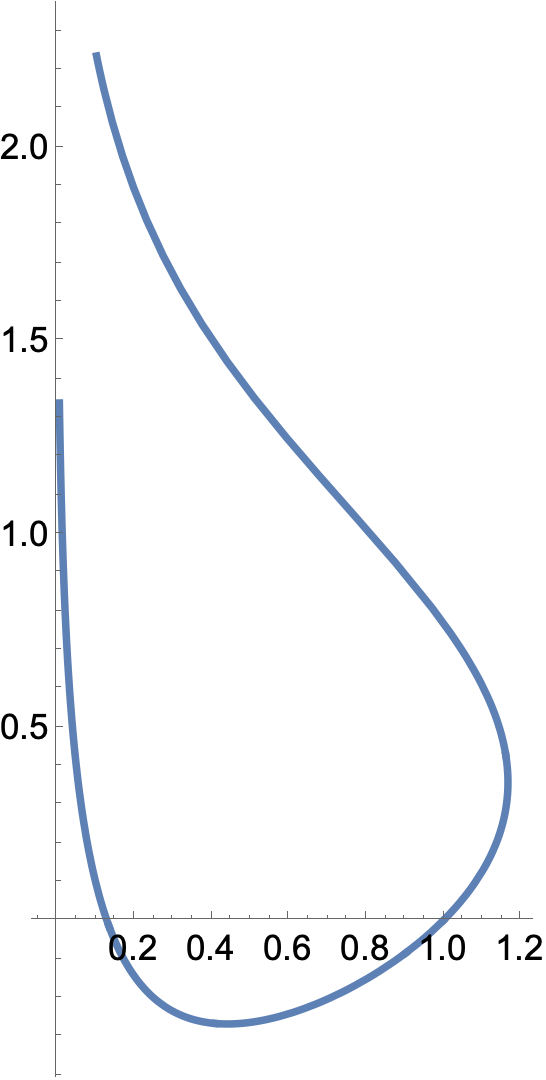}
\end{center}
\caption{The generating curve $\alpha$: $\varphi$ is constant (left)
and for the initial conditions 
$\alpha(0)=(1,0)$ and $\alpha'(0)=(0,1)$ (right)}
\label{fig:4-1}
\end{figure}

\item  From \eqref{nh-a}, we have $\langle N,\partial_\theta\rangle=\frac{x'}{\sqrt{2}\Phi}$. Then Eq. \eqref{eq2} 
\[t^2(x'\theta''-\theta'x'')=x'((x'+t\theta')^2+t^2\theta'^2).\]
This is   equivalent to
\[
t^2( x'\theta''-\theta'x'' -2x'\theta'^2)-2x'^2\theta' t-x'^3=0.\]
Because the arbitrariness of $t$, we must have $x'=0$ identically, which implies that
$\Sigma=\Sigma_{x_0}$ for a certain constant $x_0$.

\item From \eqref{nh-k}, we have $\langle N,e_3\rangle=0$. Thus  \eqref{eq2} implies that the surface 
is minimal. For a detailed study of constant mean curvature rotational surfaces in 
$\sl$ see \cite[\S4]{ko}.
 \end{enumerate}
\end{proof}



\section{Translators by the vector field $x\partial_x+y\partial_y$}
\label{s5}
Consider the Killing vector field  
 
$$V=x\partial_x+y\partial_y.$$
In terms of the basis $B$, we have  $V=\frac{1}{2y}(xe_1+y e_2+xe_3)$.

\begin{theorem} 
\begin{enumerate}
\item The only $\nn$-invariant $V$-translators are the surfaces of type 
$\Sigma_{\theta_0}$ or they are parametrized by \eqref{p-Nsurf1}, where 
$\alpha(s)= \left(c(1+\cos(\sqrt{2}(s-s_0)),s\right)$, where $c>0$ and $s_0\in\mathbb{R}$. 
\item The only $\aa$-invariant $V$-translators are the surfaces of type $\Sigma_{\theta_0}$ or of type $\Sigma_{x_0}$ with $x_0=0$.
\item Let $\Sigma$ be a rotational surface whose generating curve is parametrized   by \eqref{xy}. If $\Sigma$ is a $V$-translator, then 
 \begin{equation}\label{nv3}
 \varphi'=-\frac{y\cos\varphi+x\sin\varphi}{y}.
 \end{equation} 
\end{enumerate}
\end{theorem}

\begin{proof} 
\begin{enumerate}
\item 
Obviously, $\Sigma_{\theta_0}$ is a solution of the problem. Indeed, the normal
to $\Sigma_{\theta_0}$ is $N=\frac{1}{\sqrt{2}}(e_1-e_3)$ which is orthogonal to $V$.
On the other hand, $\Sigma_{\theta_0}$ is minimal.

Suppose that $\theta'\neq0$, hence we can choose $\theta(s)=s$. The equation \eqref{eq2}
becomes
$$
y'^2=2y(y''+y).
$$
With the change of variable $f(s)=\frac{y'(s)}{y(s)}$, this equation writes then as 
$$
2f'+f^2+2=0.
$$
The   solution of this equation is $f(s)=-\sqrt{2}\tan\left(\frac{s-s_0}{\sqrt{2}}\right)$, for
a certain constant $s_0$ obtained from the initial conditions. Then we have proved 
$$\frac{y'}{y}=-\sqrt{2}\tan\left(\frac{s-s_0}{\sqrt{2}}\right).$$
Consequently, we find
$$
y(s)=c\left(1+\cos(\sqrt{2}(s-s_0))\right),\quad c>0.
$$
  The domain of $s$ is such that 
$\frac{s-s_0}{\sqrt{2}}\in(-\frac{\pi}{2},\frac{\pi}{2})$.

\item  By using the parametrization \eqref{p-conoid}, we obtain 
$\langle N,V\rangle=-\frac{ x\theta'}{\Phi}$.
Equation \eqref{eq2} becomes 
$$
2t^2(x'\theta''-\theta'x'')=-x\theta'((x'+2t\theta')^2+x'^2).
$$
This equation writes as 
$$t^2(x'\theta''-\theta'x''+2x\theta'^3)+2txx'\theta'^2+x\theta'x'^2=0.$$
As $t$ is arbitrary, we must have either $\theta'(s)=0$ for all $s$, that is $\Sigma=\Sigma_{\theta_0}$,
  or $\theta'\neq0$ and this implies $x=0$, that is $\Sigma=\Sigma_{x_0}$ 
with $x_0=0$.

 \item Since $\langle N,V\rangle=\frac{1}{2y}(-x\sin\varphi+y\cos\varphi)$, then \eqref{eq2}  writes as in \eqref{nv3}.

Obviously, $\sin\varphi\neq0$ on a certain interval.    From the relation above we deduce
 $$
 \frac{x}{y}=-\frac{\varphi'+\cos\varphi}{\sin\varphi}, 
 \quad\textrm{for } \sin\theta\neq0.
 $$
  \end{enumerate}
\end{proof}

Let us emphasize two particular situations of Eq.  \eqref{nv3}. 

\begin{corollary} Let $\Sigma$ be a rotational $V$-translator whose generating curve is parametrized   by \eqref{xy}.
\begin{enumerate}
\item If $\varphi$ is constant, then $\Sigma$ is the surface of Ex. \ref{ex2} with $\varphi(s)=\pi/2$.
\item If $\Sigma$ has constant mean curvature $H$, then $H=0$ and the surface is of previous item (1).  
\end{enumerate}
\end{corollary}

\begin{proof}
The first part was proved in Ex. \ref{ex2}. Suppose now that the mean curvature $H$ is constant.   We must have
$$
\varphi'=2(H-\cos\varphi)\quad\textrm{and}\quad
\varphi'=-\frac{y\cos\varphi+x\sin\varphi}{y}.
$$
It follows that $2H=\cos\varphi-\frac{x}{y}\sin\varphi$.
Taking the derivative and after some manipulations we obtain
$$
(\sin\varphi+\frac{x}{y}\cos\varphi)\varphi'=
-2\sin\varphi(\cos\varphi-\frac{x}{y}\sin\varphi).
$$
Thus
$H\sin\varphi=\frac{x}{y}\left(1-H\cos\varphi\right).$ Multiply by
$\sin\varphi$ and replace $\frac{x}{y}\sin\varphi=\cos\varphi-2H$ to obtain
$$
3H=(1+2H^2)\cos\varphi.
$$ 
It follows that $H=0$ and  $\varphi$ is constant with $\varphi(s)=\frac{\pi}{2}$.
\end{proof}

\section{Translators by the vector field $\frac12(x^2-y^2)\partial_x+xy\partial_y $}\label{s6}
Consider the Killing vector field  
 
$$W=\frac12(x^2-y^2)\partial_x+xy\partial_y=\frac{1}{2y}\left(\frac12(x^2-y^2) e_1+xy e_2+\frac12(x^2-y^2)e_3\right).$$

\begin{theorem} 
\begin{enumerate}
\item The only  $\nn$-invariant $W$-translators are the surfaces  of type $\Sigma_{\theta_0}$.
\item The  only  $\aa$-invariant   $W$-translators are surfaces of type $\Sigma_{\theta_0}$.
\item Let $\Sigma$ be a rotational surface whose generating curve is parametrized by   \eqref{xy}. If $\Sigma$ is a $W$-translator, then 
\begin{equation}\label{rot6}
\varphi'=(x-2)\cos\varphi-\frac{1}{2y}(x^2-y^2)\sin\varphi.
\end{equation}
\end{enumerate}
\end{theorem}

\begin{proof} 
\begin{enumerate}
\item From \eqref{nh-n} we have 
$\langle N,W\rangle=-\frac{t\theta'}{2y\Phi}$, while $H$ depends only on $s$.
This implies $\theta'=0$ identically and thus $\theta(s)=\theta_0$ is a constant function. This proves that the surface is of type $\Sigma_{\theta_0}$.

 
\item We have
$$\langle N,W\rangle=-\frac{x^2-t^2}{2\Phi}\theta'.$$ 
Consequently, equation \eqref{eq2} writes as
$$
 4t^2(x'\theta''-x''\theta')+(x^2-t^2)\theta'\left((x')^2 (x'+ 2t\theta')^2\right)=0.
$$
The arbitrariness of $t$ implies that all coefficients appearing above vanish.
Since the coefficient of $t^4$ is zero, we must have $\theta'=0$. 
Hence $\Sigma=\Sigma_{\theta_0}$.


 \item We have 
\begin{equation}\label{nw}
 \langle N,W\rangle=\frac{1}{2y}\left(\frac12(y^2-x^2)\sin\varphi+xy\cos\varphi\right),
\end{equation}
  then \eqref{eq2} becomes
  \begin{equation}
  \label{eq:rot-W}
\frac{\varphi'}{2}+\cos\varphi=\frac{y^2-x^2}{4y}\sin\varphi+\frac{x}{2}\cos\varphi.  
  \end{equation}

 \end{enumerate}
\end{proof}
Equation \eqref{rot6}, together equations \eqref{xy} are difficult to solve. A particular case to consider is when the surface has constant mean curvature and it is natural to ask if there exist rotational $W$-translators with constant mean curvature. The answer is no as  we will prove in the   next result. This is a consequence that it is possible to find explicit parametrizations of rotational surfaces with constant mean curvature. 

\begin{corollary} There are no rotational $W$-translators with constant mean curvature.
\end{corollary}

\begin{proof} Suppose that $H$ is constant which, without loss of generality, we can assume to be non-negative. Then \eqref{nh-k2} implies $\varphi'(s)=2H-2\cos\varphi(s)$. The solution of this ODE is, up to translations in the $s$-parameter:
\begin{equation*}
\begin{split}
\varphi(s)&=-2 \tan ^{-1}\left(\frac{(1-H) \tanh \left(\sqrt{1-H^2} s\right)}{\sqrt{1-H^2}}\right),\quad (0\leq H<1),\\
\varphi(s)&=-2 \cot ^{-1}(2 s), \quad {\textrm{or}}\quad \varphi(s)=0, \quad (H=1),\\
\varphi(s)&=2 \tan ^{-1}\left(\frac{(H-1) \tan \left(\sqrt{H^2-1} s\right)}{\sqrt{H^2-1}}\right),\quad (H>1).
\end{split}
\end{equation*}
Once we have $\varphi$, we can explicitly integrate \eqref{xy}. In the following expressions,   $c$ is a constant of integration.
\begin{enumerate}
\item Case $0\leq H<1$. 
\begin{equation*}
\begin{split}
x(s)&=-c\frac{ \sqrt{1-H^2} \sinh (2 \sqrt{1-H^2} s) }{\cosh (2 \sqrt{1-H^2} s)+H},\\
y(s)&=c(H-\cos \left(2 \tan ^{-1}\left(\frac{(1-H) \tanh (\sqrt{1-H^2} s)}{\sqrt{1-H^2}}\right)\right)).\\
\end{split}
\end{equation*} 
\item Case $H=1$.  If $\varphi(s)=0$ we obtain the surface $\Sigma_{y_0}$ which is not a
 translator. For the other value of $\varphi$ we get
\begin{equation*}
\begin{split}
x(s)&=-\frac{1}{2} c \sin (2 \cot ^{-1}(2 s)),\\
y(s)&=\frac{c}{4 s^2+1}.\\
\end{split}
\end{equation*}
\item Case $H>1$. 
\begin{equation*}
\begin{split}
x(s)&=c\frac{\sqrt{H^2-1} \sin (2 \sqrt{H^2-1} s)}{\cos (2 \sqrt{H^2-1} s)+H},\\
y(s)&=c(H-\cos \left(2 \tan ^{-1}\left(\frac{(H-1) \tan (\sqrt{H^2-1} s)}{\sqrt{H^2-1}}\right)\right)).\\
\end{split}
\end{equation*}
\end{enumerate}
 Finally, we compute $\langle N,W\rangle$ using \eqref{nw} and we check that, indeed, $H\not=\langle N,W\rangle$:
 \begin{equation*}
 \begin{split}
 \langle N,W\rangle&=-\frac{1}{4} c \sqrt{1-H^2} \sinh (2 \sqrt{1-H^2} s),\quad (0\leq H<1),\\
  \langle N,W\rangle&=0,\quad (H=1),\\
   \langle N,W\rangle&=\frac{1}{4} c\sqrt{H^2-1} \sin (2 \sqrt{H^2-1} s),\quad (H>1).
   \end{split}
   \end{equation*}
\end{proof}

\section*{Acknowledgements}  

Rafael L\'opez  is a member of the IMAG and of the Research Group ``Problemas variacionales en geometr\'{\i}a'',  Junta de Andaluc\'{\i}a (FQM 325). This research has been partially supported by MINECO/MICINN/FEDER grant no. PID2020-117868GB-I00, and by the ``Mar\'{\i}a de Maeztu'' Excellence Unit IMAG, reference CEX2020-001105- M, funded by MCINN/AEI/10.13039/501100011033/ CEX2020-001105-M. Marian Ioan Munteanu is thankful to Romanian Ministry of Research, Innovation and Digitization, within Program 1 – Development of the national RD system, Subprogram 1.2 – Institutional Performance – RDI excellence funding projects, Contract no.11PFE/30.12.2021, for financial support.

 \section*{Declarations statements}

\noindent{\bf Author contributions.} All authors contributed equally in conception, design and preparation of the manuscript. 

\noindent{\bf Competing interests.} The author declares no competing interests.

\noindent {\bf Data availability.} No data sets were generated or analyzed during the current study.

\end{document}